\documentclass[11pt]{article}
\usepackage{amssymb, caption, curves, epic, mathdots}
\newtheorem{lm}{Lemma}[section]
\newtheorem{thm}{Theorem}[section]
\newtheorem{cor}{Corollary}[section]

\begin{document}
\title{\bf Graphs of order $n$ and diameter $2(n-1)/3$ minimizing the spectral radius}
\author{Jingfen Lan\\
Center for Combinatorics, LPMC-TJKLC\\
Nankai University, Tianjin 300071, China\\
{\tt jflan@sina.cn}\\
and\\
Lingsheng Shi\thanks{Project 91338102 supported by National Natural Science Foundation of China.}\\
Department of Mathematical Sciences\\
Tsinghua University, Beijing 100084, China\\
{\tt lshi@math.tsinghua.edu.cn}
}
\date{}
\maketitle

\begin{abstract}
  The spectral radius of a graph is the largest eigenvalue of its adjacency matrix. A minimizer graph is such that minimizes the spectral radius among all connected graphs on $n$ vertices with diameter $d$. The minimizer graphs are known for $d\in\{1,2\}\cup [n/2,2n/3-1]\cup\{n-k\mid k=1,2,...,8\}$. In this paper, we determine all minimizer graphs for $d=2(n-1)/3$.\\

\noindent {\em AMS classification}: 05C35; 05C50; 05E99; 94C15\\
\noindent {\em Keywords:} Diameter; Spectral radius
\end{abstract}

\section{Introduction}

All graphs considered in this paper are undirected and simple. Let $G$ be a graph. The greatest distance between any two vertices in $G$ is the {\em diameter} of  $G$, denoted by $d(G)$, or simply by $d$. An {\em internal path} of a graph is a path whose internal vertices have degree 2 and the two end vertices have degree at least 3. An internal path is closed if its two end vertices coincide. The {\em characteristic polynomial} of
$G$, simply denoted by $\phi _G$, is defined by $\phi _G(\lambda )=\det (\lambda I-A(G))$, where $A(G)$ is
the adjacency matrix of $G$ and $I$ is an identity matrix. The largest root of $\phi _G$ is the {\em spectral radius} of $G$, denoted by $\rho (G)$.

Hoffman and Smith \cite{Hoffman,HS,smith} completely determined all connected graphs $G$ with $\rho (G)\le 2$. Cvetkovi\'c et al. \cite{CDG}, Brouwer and Neumaier \cite{BN} characterized all connected graphs $G$ with $2 < \rho(G) \le\sqrt{2 + \sqrt{5}}$.

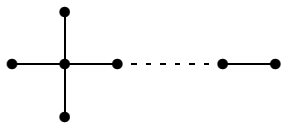
\begin{figure}[h]
\begin{center}
\setlength{\unitlength}{7mm}
\begin{picture}(5,2.2)
\multiput(0,1)(1,0){3}{\circle*{0.2}}
\multiput(1,0)(0,2){2}{\circle*{0.2}}
\multiput(4,1)(1,0){2}{\circle*{0.2}}
\put(0,1){\line(1,0){2}}
\put(1,0){\line(0,1){2}}
\put(4,1){\line(1,0){1}}
\dashline{.1}(2,1)(4,1)
\end{picture}
\end{center}
\caption{A dagger}\label{dagger}
\end{figure}

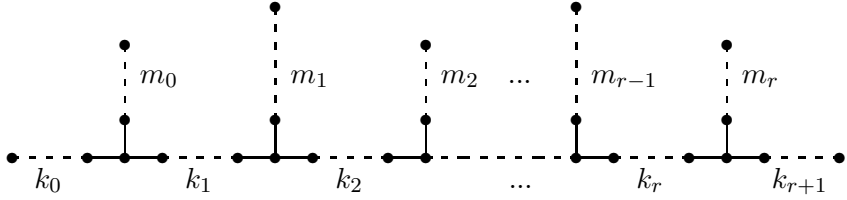
\begin{figure}[h]
\begin{center}
\setlength{\unitlength}{1cm}
\begin{picture}(11,3)
\multiput(0,.5)(1,0){6}{\circle*{.15}}
\multiput(1.5,.5)(2,0){5}{\circle*{.15}}
\multiput(8,.5)(1,0){3}{\circle*{.15}}
\put(11,.5){\circle*{.15}}
\multiput(0,.5)(2,0){6}{\dashline{.1}(0,0)(1,0)}
\multiput(5.5,.5)(1.5,0){2}{\dashline{.1}(0,0)(.5,0)}
\multiput(1,.5)(2,0){2}{\line(1,0){1}}
\multiput(5,.5)(2.5,0){2}{\line(1,0){.5}}
\put(9,.5){\line(1,0){1}}
\multiput(1.5,1)(0,1){2}{\circle*{.15}}
\multiput(3.5,1)(0,1.5){2}{\circle*{.15}}
\multiput(5.5,1)(0,1){2}{\circle*{.15}}
\multiput(7.5,1)(0,1.5){2}{\circle*{.15}}
\multiput(9.5,1)(0,1){2}{\circle*{.15}}
\multiput(1.5,.5)(2,0){5}{\line(0,1){.5}}
\multiput(1.5,1)(4,0){3}{\dashline{.1}(0,0)(0,1)}
\multiput(3.5,1)(4,0){2}{\dashline{.1}(0,0)(0,1.5)}
\put(.3,.1){$k_0$}\put(2.3,.1){$k_1$}\put(4.3,.1){$k_2$}\put(6.6,.1){...}\put(8.3,.1){$k_r$}\put(10.1,.1){$k_{r+1}$}
\put(1.7,1.5){$m_0$}\put(3.7,1.5){$m_1$}\put(5.7,1.5){$m_2$\quad ...}\put(7.7,1.5){$m_{r-1}$}\put(9.7,1.5){$m_r$}
\end{picture}
\caption{The open quipu $P_{(k_0,k_1,...,k_r,k_{r+1})}^{(m_0,m_1,...,m_r)}$}\label{opqui}
\end{center}
\end{figure}

A {\em dagger} is obtained by adding a pendent path to the center of a star of order 4, see Figure~\ref{dagger}; an {\em open quipu} is a tree with maximum degree $3$ such that all vertices of degree $3$ lie on a path; a {\em closed quipu} is a unicyclic graph with maximum degree $3$ such that all vertices of degree $3$ lie on the cycle. An open (or closed) quipu can be written in the form of $P_{(k_0,k_1,...,k_r,k_{r+1})}^{(m_0,m_1,...,m_r)}$ (or $C_{(k_1,...,k_r)}^{(m_1,...,m_r)}$) with all $k_i,m_i\ge 0$ and $r\ge 0$ (or $r\ge 1$), where for $r\ge 1$ and $1\le i\le r$, $k_i$ measures the number of internal vertices on the $i$th internal path, while $k_0$, $k_{r+1}$, $m_0$ and $m_i$ stand for the lengths of the indicated pendent paths respectively, see Figures~\ref{opqui}~and~\ref{cloqui}. These terminologies were first introduced by Woo and Neumaier \cite{WN} for the following result.
\begin{figure}[h]
\begin{center}
\setlength{\unitlength}{1cm}
\begin{picture}(7,7)
\put(3.5,3.5){\thicklines\bigcircle[4]{4}}
\multiput(3.5,0)(0,1){2}{\circle*{0.15}}
\multiput(3.5,0)(0,6){2}{\dashline{.1}(0,0)(0,1)}
\multiput(3.5,1)(0,4.5){2}{\line(0,1){.5}}
\multiput(3.5,1.5)(0,4){2}{\circle*{0.15}}
\multiput(3.5,6)(0,1){2}{\circle*{0.15}}
\multiput(1.05,1.05)(.7,.7){2}{\circle*{0.15}}
\dashline{.1}(1.05,1.05)(1.75,1.75)
\put(1.6,1.6){\line(1,1){.51}}
\put(3.9,5.45){\circle*{.15}}
\put(3.5,3.5){\thicklines\arc(.4,1.95){12}}
\dashline{.1}(5.25,5.25)(5.95,5.95)
\put(2.1,2.1){\circle*{.15}}
\put(3.5,3.5){\thicklines\arc(-1.4,-1.4){12}}
\multiput(4.9,4.9)(.35,.35){2}{\circle*{0.15}}
\put(4.9,4.9){\line(1,1){.51}}
\put(6,6){\circle*{.15}}
\multiput(5.5,3.5)(.5,0){2}{\circle*{0.15}}
\put(5.5,3.5){\line(1,0){.5}}
\dashline{.1}(6,3.5)(7,3.5)
\put(7,3.5){\circle*{.15}}
\put(4.6,5.15){\circle*{.15}}
\put(5.1,4.72){\circle*{.15}}
\put(3.5,3.5){\thicklines\arc(1.6,1.22){20}}
\put(5.48,3.9){\circle*{.15}}
\put(3.5,3.5){\thicklines\arc(2.02,0){12}}
\put(3.1,1.55){\circle*{.15}}
\put(2.4,1.8){\circle*{.15}}
\put(3.5,3.5){\thicklines\arc(-.4,-1.98){12}}
\put(2.8,6.3){$m_r$}\put(4.2,5.5){$k_1$}\put(5.2,5.8){$m_1$}\put(5.5,4.3){$k_2$}\put(6.2,3.8){$m_2$}
\put(3.7,.4){$m_{i-1}$}\put(1.4,1){$m_i$}\put(2.4,1.2){$k_i$}
\put(.8,3){$\vdots $}\put(1.2,5){$\iddots $}\put(5.2,1.2){$\iddots $}
\end{picture}
\end{center}
\caption{The closed quipu $C_{(k_1,...,k_r)}^{(m_1,...,m_r)}$}\label{cloqui}
\end{figure}
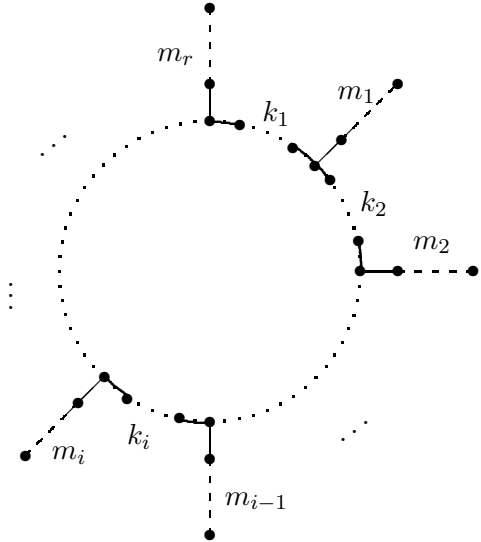
\begin{lm}{\rm \cite{WN}}\label{WandN}
 A graph $G$ whose spectral radius satisfies $2<\rho(G)\le3/\sqrt{2}$ is either an open quipu, a closed quipu, or a dagger.
\end{lm}

A {\em minimizer} graph of order $n$ with diameter $d$ is such a graph that has the minimal spectral radius among all simple connected graphs on $n$ vertices with diameter $d$. The problem to determine the minimizer graphs was raised by van Dam and Kooij \cite{DK} concerning a model of virus propagation in networks. They solved this problem explicitly for $d\in\{1, 2, \lfloor n/2\rfloor, n-3, n-2, n-1\}$, where the minimizer graph is a complete graph for $d=1$, a star for $d=2$ and $n$ large enough, a cycle for $d=\lfloor n/2\rfloor $ and $n>6$, a path for $d=n-1$, $P_{(1,n-3)}^{(1)}$ for $d=n-2$, and $P_{(1,n-6,1)}^{(1,1)}$ for $d=n-3$. All the minimizer graphs on $n\le 20$ vertices were also obtained in \cite{DK} via computer aid. The minimizer graphs are not unique in general. Later, they were further determined for $d=n-4$ by Yuan-Shao-Liu \cite{YSL}, for $d=n-5$ by Cioab\v{a}-van~Dam-Koolen-Lee \cite{CDK}, and for $d=n-6, n-7, n-8$ by Lan-Lu-Shi \cite{LLS}, which turn out all to be open quipus. Recently, Lan and Lu \cite{LL} determined the minimizer graphs for $n/2\le d\le (2n-3)/3$, where the closed quipu $C_{(n-d-1,n-d-1)}^{(d-\lfloor\frac{n}{2}\rfloor,d-\lceil\frac{n}{2}\rceil)}$ is the unique minimizer graph for $n/2\le d\le (2n-5)/3$, the closed quipus $C_{(n-d-1,n-d-1)}^{(i,2d-n-i)}$ for all integers $0\le i\le 2d-n$ form the set of all the minimizer graphs for $d=(2n-4)/3$, and there are exactly two minimizer graphs for $d=(2n-3)/3$: the closed quipu $C^{(k-2)}_{(2k+1)}$ and the open quipu $P_{(1,k-2,k-1)}^{(1,k-1)}$ with $k=n/3$.

In this paper, we determine the minimizer graphs for $d=2(n-1)/3$ which form a family of open quipus. This confirms a conjecture in \cite{LL}. Symmetric to the case $d=(2n-4)/3$, our result indicates that the minimizer graphs have a phase transition occurring at $d=(2n-3)/3$ from closed quipus to open quipus.

\begin{thm}\label{mainthm} For $k\in\mathbb N$, the minimizer graphs of order $3k+1$ and diameter $2k$ are exactly $P^{(i,j)}_{(i,i+j-1,j)}$ for all $0\le i\le j$ satisfying $i+j=k$.
\end{thm}

The proof of Theorem~\ref{mainthm} mainly relies on Lemma~\ref{WandN}. The graphs stated in Theorem~\ref{mainthm} are known (by Lemma~\ref{mm'}) with spectral radius less than $3/\sqrt 2$. Then three more steps determine the minimizer graphs.

Step 1. Exclude the family of daggers and closed quipus for minimizer graphs by the diameter condition for graphs with $2<\rho(G)<3/\sqrt 2$ (mainly by Lemma~\ref{SpD}).

Step 2. Exclude most open quipus by refining the lengths of internal paths and pendent paths (by Lemmas~\ref{sufc} and~\ref{sufc2}).

Step 3. Deal with several open quipus left to be excluded. For this purpose, some lemmas shown in the next section will be used. The proof is given in the last section.

\section{Preliminaries}

For any vertex $v$ in a graph $G$, let $N(v)$ be the neighborhood of $v$. Let $G-v$ be the remaining graph of $G$ after deleting the vertex $v$ (and all edges incident to $v$). Similarly, $G-u-v$ is the remaining graph of $G$ after deleting the two vertices $u$ and $v$.

\begin{lm}{\rm\cite{Sch}}\label{lm2.1}
 Let $G$ be a graph, $v$ be a vertex of $G$, and $e=uv$ be an edge of $G$. Suppose $v$ and $e$ are not contained in any cycle of $G$. Then the characteristic polynomial $\phi _G$ satisfies
\begin{eqnarray*}
\phi _G&=&\lambda\phi _{G-v}-\sum_{w\in N(v)}\phi _{G-w-v},\\
\phi _G&=&\phi _{G-e}-\phi _{G-u-v}.
\end{eqnarray*}
\end{lm}
The facts stated in the following lemma can be used to compare the spectral radii of two graphs.
\begin{lm}\label{lm2.2} Let $G_1$ and $G_2$ be two
graphs. Then the following statements hold.
\begin{enumerate}
\item {\rm\cite{CDS}} If $G_2$ is a proper subgraph of $G_1$, then $\rho(G_1)>\rho(G_2)$.
\item {\rm\cite{LF}} If $G_1$ is connected and $G_2$ is a proper spanning subgraph of $G_1$, then  $\rho(G_1)>\rho(G_2)$ and $\phi_{G_2}(\lambda)>\phi_{G_1}(\lambda)$ for all $\lambda\ge\rho(G_1)$.
\item If $\phi_{G_2}(\lambda)>\phi_{G_1}(\lambda)$ for all $\lambda\ge{\rho(G_1)}$, then $\rho(G_2)<\rho(G_1)$.
\item If $\phi_{G_1}(\rho(G_2))<0$, then $\rho(G_1)>\rho(G_2)$.
\end{enumerate}
\end{lm}

\begin{lm}{\rm\cite{HS}}\label{lm2.4}
Let $uv$ be an edge of a connected graph $G$ of order $n$, and denote by $G_{u,\,v}$ the graph of order $n+1$ obtained from $G$ by subdividing the edge $uv$ once, i.e., replacing the edge $uv$ by a new vertex $w$ and two new edges $uw, vw$. Then the following two properties hold.
\begin{enumerate}
\item[\rm (i)] If $uv$ does not belong to an internal path of $G$ and $G\ne
C_n$, then $\rho(G_{u,\,v})>\rho(G)$.
\item[\rm (ii)] If $uv$ belongs to an internal path of $G$ and $G\ne P_{(1,n-6,1)}^{(1,1)}$, then $\rho(G_{u,\,v})<\rho(G)$.
\end{enumerate}
\end{lm}

The following lemma indicates the effect of edge transfers on the spectral radii of graphs. The result for $k-l\ge j$ were stated in \cite{LF} without a proof. For completeness, we include a proof here.
\begin{lm}{\rm\cite{LF}}\label{trans}
Let $j,k\ge 0$ and $l>0$ be integers. Let $u$ and $v$ be two vertices (possibly $u=v$ for $j=0$) of degree at least 2 and connected by an induced path of length $j$ in a graph $G$. Denote by $G^{(j)}_{k,l}$ the graph obtained from $G$ by adding two pendent paths of lengths $k$ and $l$ to vertices $u$ and $v$ respectively, see Figure~\ref{Gjkl}. If $k-l\ge j-1$, then
\begin{eqnarray}
\phi_{G^{(j)}_{k,l}}(\lambda )&\le &\phi_{G^{(j)}_{k+1,l-1}}(\lambda )\mbox{ for }\lambda\ge\rho\left (G^{(j)}_{k+1,l-1}\right ),\label{p}\\
\rho\left (G^{(j)}_{k,l}\right )&\ge &\rho\left (G^{(j)}_{k+1,l-1}\right ),\label{r}
\end{eqnarray}
with each equality if and only if $j=0$ and $k=l-1$.
\end{lm}
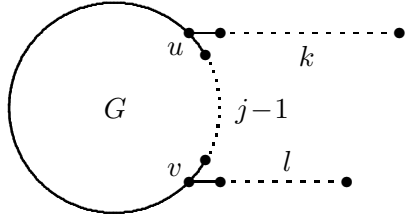
\begin{figure}[h]
\begin{center}
\setlength{\unitlength}{7mm}
\begin{picture}(8,4)
\put(2,2){\thicklines\arc(1.732,1){300}}
\put(2,2){\thicklines\arc[4](1.732,-1){60}}
\put(3.732,3){\circle*{0.2}}
\put(3.732,1){\circle*{0.2}}
\put(3.414,3.414){\circle*{0.2}}
\put(3.414,.596){\circle*{0.2}}
\put(3.414,3.414){\dashline{.1}(0,0)(4,0)}
\put(3.414,.596){\dashline{.1}(0,0)(3,0)}
\put(3.414,3.414){\line(1,0){.6}}
\put(3.414,.596){\line(1,0){.6}}
\put(4.014,3.414){\circle*{0.2}}
\put(4.014,.596){\circle*{0.2}}
\put(7.414,3.414){\circle*{0.2}}
\put(6.414,.596){\circle*{0.2}}
\put(4.3,1.8){$j\!-\!1$}\put(5.5,2.8){$k$}\put(5.2,.8){$l$}
\put(1.8,1.8){$G$}\put(3,3){$u$}\put(3,.7){$v$}
\end{picture}
\caption{The graph $G^{(j)}_{k,l}$}\label{Gjkl}
\end{center}
\end{figure}
\noindent{\bf Proof}. The strict inequalities in Eqs.~(\ref{p})~and~(\ref{r}) for $j=0$ and $k-l\ge 0$ were proven in \cite[Theorem 5]{LF}. It is easy to see that the equalities hold if $j=0$ and $k=l-1$ since the two graphs $G^{(0)}_{l-1,l}$ and $G^{(0)}_{l,l-1}$ are isomorphic. Now we assume $j>0$ and then $u\ne v$. Applying Lemma~\ref{lm2.1}, we have
\begin{eqnarray*}
\phi_{G^{(j)}_{k,l}}-\phi_{G^{(j)}_{k+1,l-1}} &=& \lambda\phi_{G^{(j)}_{k,l-1}}-\phi_{G^{(j)}_{k,l-2}}-\left (\lambda\phi_{G^{(j)}_{k,l-1}}-\phi_{G^{(j)}_{k-1,l-1}}\right )\\
&=& \phi_{G^{(j)}_{k-1,l-1}}-\phi_{G^{(j)}_{k,l-2}}=...\\
&=& \phi_{G^{(j)}_{k-l+1,1}}-\phi_{G^{(j)}_{k-l+2,0}}\\
&=& \phi_{G^{(j)}_{k-l,0}}-\phi_{G^{(j)}_{k-l+1,0}-v}.
\end{eqnarray*}

The graph $G^{(j)}_{k-l+1,0}-v$ has two pendent paths of lengths $k-l+1$ and $j-1$ at the vertex $u$. Deleting these two pendent paths results in a subgraph $H$ and $G^{(j)}_{k-l+1,0}-v=H^{(0)}_{k-l+1,j-1}$. Suppose that $k-l\ge j-1$. Applying this lemma to $H^{(0)}_{k-l,j}$, we have
\begin{eqnarray*}
\phi _{H^{(0)}_{k-l,j}}(\lambda )&\le &\phi _{H^{(0)}_{k-l+1,j-1}}(\lambda )=\phi_{G^{(j)}_{k-l+1,0}-v}(\lambda )\mbox{ for }\lambda\ge\rho\left (G^{(j)}_{k-l+1,0}-v\right ),\\
\rho\left (H^{(0)}_{k-l,j}\right )&\ge &\rho\left (G^{(j)}_{k-l+1,0}-v\right ).
\end{eqnarray*}
Note that the graph $H^{(0)}_{k-l,j}$ is isomorphic to a proper spanning subgraph of $G^{(j)}_{k-l,0}$. By Lemma~\ref{lm2.2}, we get
$$
\phi_{G^{(j)}_{k,l}}(\lambda )-\phi_{G^{(j)}_{k+1,l-1}}(\lambda )\le \phi_{G^{(j)}_{k-l,0}}(\lambda )-\phi_{H^{(0)}_{k-l,j}}(\lambda )<0,
$$
for $\lambda\ge\rho\left (G^{(j)}_{k-l,0}\right )$, which implies that $\rho\left (G^{(j)}_{k,l}\right )>\rho\left (G^{(j)}_{k+1,l-1}\right )$.
$\hfill\Box $
\medskip

Woo and Neumaier \cite{WN} noted that no (finite) graph has spectral radius  exactly $3/\sqrt{2}$ since this is not an algebraic integer. A dagger on $n$ vertices have diameter $n-3$ and its spectral radius approaches increasingly to $3/\sqrt 2$ as $n$ goes to infinity. However, some quipus have spectral radii greater than $3/\sqrt 2$. Lemma~\ref{WandN} was refined in \cite{LL} as follows.

\begin{lm}{\rm\cite{LL}}\label{SpD}
Let $G$ be a graph  on $n$ vertices $(n\ge 13)$ with spectral radius less than $3/\sqrt 2$. If $G$ is an open quipu then its diameter $d$ satisfies $d\ge (2n-4)/3$, where the bound is tight. If $G$ is a closed quipu then its diameter $d$ satisfies $n/3<d\le 2(n-1)/3$, where the lower bound is asymptotically tight and only the closed quipu $C_{(2k+3)}^{(k)}$ with $d=2k+2$ and $n=3k+4$ takes the upper bound.
\end{lm}

Let $\delta_1$ be the indicator function of being $1$, i.e., $\delta_1(x)=1$ if $x=1$ and $0$ otherwise. The following lemmas from \cite{LL} give necessary conditions for open quipus with spectral radius at most $3/\sqrt 2$.

\begin{lm} \label{sufc}
Suppose an open quipu $P_{(m_0,k_1,...,k_r,m_{r})}^{(m_0,...,m_r)}$ (with $r\ge 2$) has spectral radius less than $3/\sqrt 2$. Then the following statements hold.
\begin{enumerate}
\item For $2\le i\le r-1$, we have $k_i\geq m_{i-1}+m_i+1 -
\left\lceil \frac{\delta_1(m_{i-1})+\delta_1(m_i)}{2}\right\rceil$.
\item  We have  $k_1 \geq m_{0}+m_1-
\left\lceil \frac{3\delta_1(m_0)+\delta_1(m_1)}{2} \right\rceil
-\left\lfloor\frac{\delta_1(m_0-1)+\delta_1(m_1-1)}{2}\right\rfloor$.
\item  We have  $k_r \geq m_{r}+m_{r-1}-
\left\lceil \frac{3\delta_1(m_r)+\delta_1(m_{r-1})}{2} \right\rceil
-\left\lfloor\frac{\delta_1(m_r-1)+\delta_1(m_{r-1}-1)}{2}\right\rfloor$.
\end{enumerate}
\end{lm}
\begin{lm}\label{sufc2} Suppose that an open quipu
$P_{(m_0,k_1,...,k_r,m_{r})}^{(m_0,...,m_r)}$  (with $r\ge 2$) satisfies
\begin{enumerate}
\item $k_i\le m_{i-1}+m_i+2 - \left\lceil
\frac{\delta_1(m_{i-1})+\delta(m_i)}{2}
\right\rceil $ for $2\le i\le r-1$;
\item $k_1\leq m_{0}+m_1
-\left\lceil \frac{3\delta_1( m_{0}) +
\delta_1( m_{1}) +\delta_1( m_{0}-1) }{2}\right\rceil
$.
\item $k_{r}\le m_{r-1}+m_r
-\left\lceil \frac{3\delta_1( m_{r}) +
\delta_1( m_{r-1}) +\delta_1( m_{r}-1) }{2}\right\rceil
$.
\end{enumerate}
Then we have $\rho\left (P_{(m_0,k_1,...,k_r,m_r)}^{(m_0,...,m_r)}\right )>3/\sqrt 2$.
\end{lm}

Denote by $\rho _{k}$ the spectral radius of $P^{(k)}_{(k,k)}$. Then $\rho _1=\sqrt 3$ and $\rho _2=2$. Note that $P^{(k)}_{(k,k)}$ is a proper subgraph of $P^{(k+1)}_{(k+1,k+1)}$. By Lemma~\ref{lm2.2} and \cite[Lemma 3]{WN}, we have
\begin{equation}\label{rho}
\rho _k<\rho _{k+1}<3/\sqrt 2.
\end{equation}
Moreover, the following lemma from \cite{LL} shows that the graphs we desire in Theorem~\ref{mainthm} share the same spectral radius.
\begin{lm}\label{mm'}{\rm \cite{LL}}
For any non-negative integers $i, j$ satisfying $i+j\ge 2$, all open quipus
$P^{(i,j)}_{(i,i+j-1,j)}$ and all closed quipus $C^{(i-1,j-1)}_{(i+j+1,i+j+1)}$ have the same spectral radius $\rho _{i+j}$. \end{lm}

Let $v$ be a vertex of graph $G$. In \cite{LLS}, a {\em rooted} graph $(G, v)$ was defined as the graph $G$ together with the designated vertex $v$ as a {\em root}, and we introduced two parameters $p_{(G,v)}$ and $q_{(G,v)}$ satisfying
\begin{eqnarray*}
 \phi_G&=&p_{(G,v)}+q_{(G,v)},\\
 \phi_{G-v}&=&x_2p_{(G,v)}+x_1q_{(G,v)}.
\end{eqnarray*}
Here $x_1$ and $x_2$ are the two roots of the equation $x^2-\lambda x +1=0$, namely
$$x_1=\frac{\lambda-\sqrt{\lambda^2-4}}{2}~\mbox{ and }~x_2=\frac{\lambda+\sqrt{\lambda^2-4}}{2}.$$
The fact $x_1+x_2=\lambda$, $x_1x_2=1$ will be used deliberately. In this paper, we always assume $\lambda> 2$, then $x_1< 1< x_2$. Thus $p_{(G,v)}$ and $q_{(G,v)}$ are well defined, and
$${p_{(G,v)}\choose q_{(G,v)}}=\frac{1}{x_2-x_1}\left (\begin{array}{rr} -x_1 & 1\\ x_2 & -1 \end{array}\right ){\phi_G\choose \phi_{G-v}}.$$
Let $P_n$ denote a path of order $n$. As an example in \cite[Section 2.2]{LL}, we have
\begin{equation}\label{eq:p2k+1}
{p_{(P_{2k+1},*)}\choose q_{(P_{2k+1},*)}}=
\frac{x_2^{k+1}-x_1^{k+1}}{(x_2-x_1)^3}
{x_2^{k-1}-2x_1^{k+1}+x_1^{k+3}\choose x_1^{k-1}-2x_2^{k+1}+x_2^{k+3}}.
\end{equation}
where $*$ stands for the center of the odd path $P_{2k+1}$ for $k\ge 0$.

\begin{figure}[h]
\begin{center}
\setlength{\unitlength}{7mm}
\begin{picture}(9,2)
\multiput(0,0)(1,0){4}{\circle*{0.2}}
\multiput(1,1)(2,0){2}{\circle*{0.2}}
\put(0,0){\line(1,0){3}}
\multiput(1,0)(2,0){2}{\line(0,1){1}}
\put(3.2,-.1){$u$}
\multiput(6,0)(1,0){4}{\circle*{0.2}}
\multiput(8,1)(0,1){2}{\circle*{0.2}}
\put(6,0){\line(1,0){3}}
\put(8,0){\line(0,1){2}}
\put(9.2,-.1){$v$}
\end{picture}
\end{center}
\caption{Two rooted graphs $\left (P_{(1,3)}^{(1)},u\right )$ and $\left (P_{(2,1)}^{(2)},v\right )$}\label{two}
\end{figure}
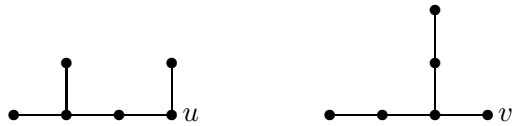

Let $t_{(G,v)}:=q_{(G,v)}/p_{(G,v)}$. It was shown in \cite{LLS} that $t_{(G,v)}$ plays an important role on the spectral radii of open quipus.
\begin{lm}\label{compare} Let $u$ and $v$ be the roots of $P_{(1,3)}^{(1)}$ and $P_{(2,1)}^{(2)}$ respectively as shown in Figure~\ref{two}. Then we have
\begin{equation}\label{tt1}
t_{\left (P_{(1,3)}^{(1)},u\right )}(\lambda )<t_{\left (P_{(2,1)}^{(2)},v\right )}(\lambda )\mbox{ for }\lambda >2.
\end{equation}
\end{lm}

\noindent{\bf Proof.} By \cite[Lemma 2.6]{LLS} and Eq.~(\ref{eq:p2k+1}), we have
\begin{eqnarray*}
\left(\begin{array}{c} p_{\left (P_{(1,3)}^{(1)},u\right )}\\ q_{\left (P_{(1,3)}^{(1)},u\right  )}\end{array}\right)
&=&\frac{1}{x_2-x_1}\left (\begin{array}{ll} \lambda-x_1^3 & x_1\\ -x_2 & x_2^3-\lambda \end{array}\right )\left (\begin{array}{ll} x_1 & 0\\ 0 & x_2\end{array}\right ){p_{(P_3,*)}\choose q_{(P_3,*)}}\\
&=& \frac{\lambda}{x_2-x_1}{(\lambda-x_1^3)x_1^3+x_2^2\choose (x_2^3-\lambda)x_2^3-x_1^2},
\end{eqnarray*}
and
$$\left(\begin{array}{c} p_{\left (P_{(2,1)}^{(2)},v\right )}\\ q_{\left (P_{(2,1)}^{(2)},v\right )}\end{array}\right)=\left (\begin{array}{ll} x_1 & 0\\ 0 & x_2\end{array}\right ){p_{(P_5,*)}\choose q_{(P_5,*)}}=\frac{\lambda^2-1}{x_2-x_1}{(\lambda-x_1^3)x_1^2\choose (x_2^3-\lambda)x_2^2}.$$
Thus we obtain
\begin{eqnarray*}
t_{\left (P_{(1,3)}^{(1)},u\right )}(\lambda )&=&\frac{(x_2^3-\lambda)x_2^3-x_1^2}{(\lambda-x_1^3)x_1^3+x_2^2},\\
t_{\left (P_{(2,1)}^{(2)},v\right )}(\lambda )&=&\frac{(x_2^3-\lambda)x_2^2}{(\lambda-x_1^3)x_1^2}.
\end{eqnarray*}
It follows that Eq.~(\ref{tt1}) is equivalent to
$$\frac{(x_2^3-\lambda)x_2^3-x_1^2}{(\lambda-x_1^3)x_1^3+x_2^2}<\frac{(x_2^3-\lambda)x_2^2}{(\lambda-x_1^3)x_1^2},$$
namely,
$$(x_2^3-\lambda)(\lambda-x_1^3)(x_2-x_1)<(x_2^3-\lambda)x_2^4+(\lambda-x_1^3)x_1^4,$$
which holds by the following easy calculation,
\begin{eqnarray*}
(x_2^3-\lambda)(\lambda-x_1^3)(x_2-x_1)&=&\left [\lambda(x_2^3+x_1^3)-\lambda^2-1\right ](x_2-x_1)\\
&<&\left [(x_2^3+x_1^3)^2-\lambda^2-1\right ](x_2-x_1)\\
&=&(x_1^6+x_2^6-x_1^2-x_2^2-1)(x_2-x_1)\\
&=&x_2^7-x_2^5-x_2^3+x_1^3+x_1^5-x_1^7\\
&=&(x_2^3-\lambda)x_2^4+(\lambda-x_1^3)x_1^4,
\end{eqnarray*}
where the inequality holds by $$x_1^3+x_2^3=(x_1+x_2)\left (x_1^2-x_1x_2+x_2^2\right )=\lambda\left  [(x_1-x_2)^2+1\right ]>\lambda >2.$$
The proof is complete. $\hfill\Box$
\medskip

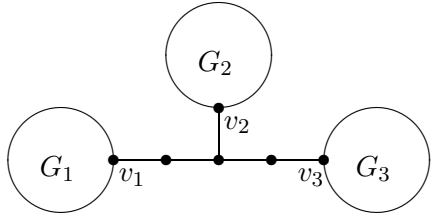
\begin{figure}[h]
\begin{center}
\setlength{\unitlength}{7mm}
\begin{picture}(8,4)
\put(4,3){\circle{2.5}}
\multiput(1,1)(6,0){2}{\circle{2.5}}
\put(4,2){\circle*{0.2}}
\multiput(2,1)(1,0){5}{\circle*{0.2}}
\put(4,1){\line(0,1){1}}
\put(2,1){\line(1,0){4}}
\put(3.6,2.7){$G_2$}
\put(0.6,.7){$G_1$}\put(6.6,.7){$G_3$}
\put(4.1,1.6){$v_2$}\put(2.1,.6){$v_1$}\put(5.5,.6){$v_3$}
\end{picture}
\end{center}
\caption{The graph $T_{G_1,G_3}^{G_2}$}\label{p3}
\end{figure}

For a vertex $v$ of graph $G$, denote by $(G,v,i)$ ($i\ge 0$) the graph obtained from $G$ by adding a pendent path of length $i$ to $v$. It is clear that $(G,v)$ can be regarded as $(G,v,0)$. Let $u$ be the other end of the pendent path in $(G,v,i)$, then by \cite[Lemma 2.6 (1)]{LLS},
\begin{eqnarray*}
\phi_{(G,v,i)} &=& (1,1){p_{[(G,v,i),u]}\choose q_{[(G,v,i),u]}}=(1,1)\left(
  \begin{array}{ll}
    x_1 & 0\\
   0 & x_2
  \end{array}
\right)^i{p_{(G,v)}\choose q_{(G,v)}}=x_1^ip_{(G,v)}+x_2^iq_{(G,v)}.
\end{eqnarray*}
Let $\alpha_{(G,v,i)}:=\phi_{(G,v,i+1)}/\phi_{(G,v,i)}$, then the following equality holds accordingly,
\begin{equation}\label{alphaequ2}
\alpha_{(G,v,i)}=\frac{\phi_{(G,v,i+1)}}{\phi_{(G,v,i)}}=\frac{x_1^{i+1}p_{(G,v)}
+x_2^{i+1}q_{(G,v)}}{x_1^ip_{(G,v)}+x_2^iq_{(G,v)}}=\frac{x_1^{2i+1}+x_2t_{(G,v)}}{x_1^{2i}+t_{(G,v)}}.
\end{equation}
For convenience, we write $\alpha_{(G,v)}$ for $\alpha_{(G,v,0)}$. Let $(G_i,v_i)$ be a (possibly empty) rooted graph for $i=1,2,3$, and let $T_{G_1,G_3}^{G_2}$ be the graph shown in Figure \ref{p3}. We have the following lemma, which indicates that the spectral radius of $T_{G_1,G_3}^{G_2}$  decreases as $\alpha_{(G_i,v_i)}$ (also $t_{(G_i,v_i)}$) increases for $i=1,2,3$.
\begin{lm}\label{trg}
The spectral radius $\rho\left (T_{G_1,G_3}^{G_2}\right )$ is the largest root of the equation
\begin{equation}\label{trgequ}
\alpha_{(G_2,v_2)}=\frac{1}{\alpha_{(G_1,v_1)}}+\frac{1}{\alpha_{(G_3,v_3)}}.
\end{equation}
Moreover, let $(G'_i,v'_i)$ be a rooted graph for $i=1,2$, then the following holds.
\begin{itemize}
\item If $\alpha _{(G_1,v_1)}\left (\rho\left (T_{G_1,G_3}^{G_2}\right )\right )>\alpha _{(G'_1,v'_1)}\left (\rho\left (T_{G_1,G_3}^{G_2}\right )\right )$, then $\rho\left (T_{G_1,G_3}^{G_2}\right )<\rho\left (T_{G'_1,G_3}^{G_2}\right )$.
\item If $\alpha _{(G_2,v_2)}\left (\rho\left (T_{G_1,G_3}^{G_2}\right )\right )>\alpha _{(G'_2,v'_2)}\left (\rho\left (T_{G_1,G_3}^{G_2}\right )\right )$, then $\rho\left (T_{G_1,G_3}^{G_2}\right )<\rho\left (T_{G_1,G_3}^{G'_2}\right )$.
\end{itemize}
\end{lm}
Lemma~\ref{trg} readily implies the following result.
\begin{cor}\label{trgequ2}{\rm \cite{Sun}}
For any pair of graphs $G_1$ and $G_2$, $\rho\Big(T_{G_1,G_1}^{G_2}\Big)=\rho\Big(T_{G_2,G_2}^{G_1}\Big)$.
\end{cor}

\noindent{\bf Proof of Lemma~\ref{trg}.} By Lemma~\ref{lm2.1}, we have
\begin{eqnarray*}
\phi_{T_{G_1,G_3}^{G_2}}
&=&\phi_{(G_1,v_1,1)}\phi_{(G_2,v_2,1)}\phi_{(G_3,v_3,1)}-\phi_{(G_1,v_1,1)}\phi_{G_2}\phi_{G_3} -\phi_{G_1}\phi_{G_2}\phi_{(G_3,v_3,1)}\\
&=&\phi_{G_2}\phi_{(G_1,v_1,1)}\phi_{(G_3,v_3,1)}\left ( \alpha_{(G_2,v_2)}-\frac{1}{\alpha_{(G_1,v_1)}}-\frac{1}{\alpha_{(G_3,v_3)}}\right ).
\end{eqnarray*}
The graphs $G_2$, $(G_1,v_1,1)$, and $(G_3,v_3,1)$ have spectral radius all  less than $\rho\left (T_{G_1,G_3}^{G_2}\right )$ since they are proper subgraphs of $T_{G_1,G_3}^{G_2}$. Thus, $\rho\left (T_{G_1,G_3}^{G_2}\right )$ must be the largest root of Eq.~(\ref{trgequ}). The rest of the lemma follows easily from Lemma~\ref{lm2.2}. $\hfill\Box $

\section{Proof of main theorem}

\noindent{\bf Proof of Theorem~\ref{mainthm}.} The theorem holds for $n:=3k+1\le 20$ as checked in \cite{DK}. So we can assume that $n>20$ and $k>6$. Lemma~\ref{mm'} together with Eq.~(\ref{rho}) implies that all graphs stated in the theorem have the same spectral radius $\rho _k\in\left (2,3/\sqrt 2\right )$.

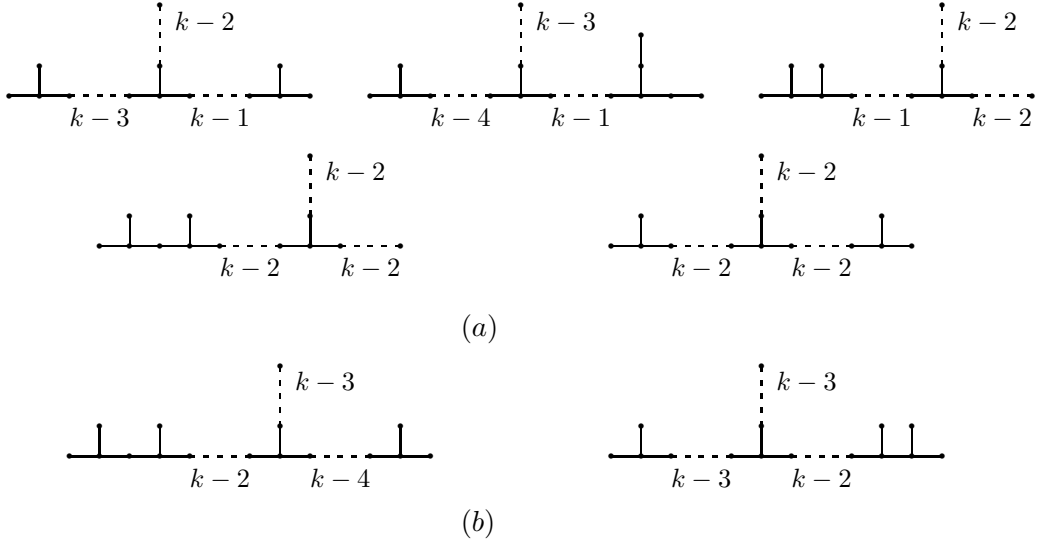
\begin{figure}[h]
\begin{center}
\setlength{\unitlength}{.4cm}
\begin{picture}(34,17)
\multiput(0,14)(1,0){3}{\circle*{0.2}}
\multiput(4,14)(1,0){3}{\circle*{0.2}}
\multiput(8,14)(1,0){3}{\circle*{0.2}}
\multiput(1,15)(4,0){3}{\circle*{0.2}}
\put(5,17){\circle*{0.2}}
\multiput(0,14)(4,0){3}{\line(1,0){2}}
\multiput(1,14)(4,0){3}{\line(0,1){1}}
\dashline{0.2}(2,14)(4,14)
\dashline{0.2}(6,14)(8,14)
\dashline{0.2}(5,15)(5,17)
\put(2,13){\small{$k-3$}}
\put(6,13){\small{$k-1$}}
\put(5.5,16.2){\small{$k-2$}}
\multiput(12,14)(1,0){3}{\circle*{0.2}}
\multiput(16,14)(1,0){3}{\circle*{0.2}}
\multiput(20,14)(1,0){4}{\circle*{0.2}}
\multiput(13,15)(4,0){3}{\circle*{0.2}}
\put(17,17){\circle*{0.2}}\put(21,16){\circle*{0.2}}
\multiput(12,14)(4,0){2}{\line(1,0){2}}
\multiput(13,14)(4,0){3}{\line(0,1){1}}
\put(20,14){\line(1,0){3}}
\put(21,14){\line(0,1){2}}
\dashline{0.2}(14,14)(16,14)
\dashline{0.2}(18,14)(20,14)
\dashline{0.2}(17,15)(17,17)
\put(14,13){\small{$k-4$}}
\put(18,13){\small{$k-1$}}
\put(17.5,16.2){\small{$k-3$}}
\multiput(25,14)(1,0){4}{\circle*{0.2}}
\multiput(30,14)(1,0){3}{\circle*{0.2}}
\multiput(26,15)(1,0){2}{\circle*{0.2}}
\multiput(31,15)(0,2){2}{\circle*{0.2}}
\put(34,14){\circle*{0.2}}
\put(25,14){\line(1,0){3}}
\put(30,14){\line(1,0){2}}
\put(31,14){\line(0,1){1}}
\multiput(26,14)(1,0){2}{\line(0,1){1}}
\dashline{0.2}(28,14)(30,14)
\dashline{0.2}(32,14)(34,14)
\dashline{0.2}(31,15)(31,17)
\put(28,13){\small{$k-1$}}
\put(32,13){\small{$k-2$}}
\put(31.5,16.2){\small{$k-2$}}
\multiput(3,9)(1,0){5}{\circle*{0.2}}
\multiput(9,9)(1,0){3}{\circle*{0.2}}
\multiput(4,10)(2,0){2}{\circle*{0.2}}
\multiput(10,10)(0,2){2}{\circle*{0.2}}
\put(13,9){\circle*{0.2}}
\put(3,9){\line(1,0){4}}
\put(9,9){\line(1,0){2}}
\put(10,9){\line(0,1){1}}
\multiput(4,9)(2,0){2}{\line(0,1){1}}
\dashline{0.2}(7,9)(9,9)
\dashline{0.2}(11,9)(13,9)
\dashline{0.2}(10,10)(10,12)
\put(7,8){\small{$k-2$}}
\put(11,8){\small{$k-2$}}
\put(10.5,11.2){\small{$k-2$}}
\multiput(20,9)(1,0){3}{\circle*{0.2}}
\multiput(24,9)(1,0){3}{\circle*{0.2}}
\multiput(28,9)(1,0){3}{\circle*{0.2}}
\multiput(21,10)(4,0){3}{\circle*{0.2}}
\put(25,12){\circle*{0.2}}
\multiput(20,9)(4,0){3}{\line(1,0){2}}
\multiput(21,9)(4,0){3}{\line(0,1){1}}
\dashline{0.2}(22,9)(24,9)
\dashline{0.2}(26,9)(28,9)
\dashline{0.2}(25,10)(25,12)
\put(22,8){\small{$k-2$}}
\put(26,8){\small{$k-2$}}
\put(25.5,11.2){\small{$k-2$}}
\put(15,6){$(a)$}
\multiput(2,2)(1,0){5}{\circle*{0.2}}
\multiput(8,2)(1,0){3}{\circle*{0.2}}
\multiput(12,2)(1,0){3}{\circle*{0.2}}
\multiput(3,3)(2,0){2}{\circle*{0.2}}
\multiput(9,3)(0,2){2}{\circle*{0.2}}
\put(13,3){\circle*{0.2}}
\put(2,2){\line(1,0){4}}
\multiput(8,2)(4,0){2}{\line(1,0){2}}
\multiput(3,2)(2,0){2}{\line(0,1){1}}
\multiput(9,2)(4,0){2}{\line(0,1){1}}
\dashline{0.2}(6,2)(8,2)
\dashline{0.2}(10,2)(12,2)
\dashline{0.2}(9,3)(9,5)
\put(6,1){\small{$k-2$}}
\put(10,1){\small{$k-4$}}
\put(9.5,4.2){\small{$k-3$}}
\multiput(20,2)(1,0){3}{\circle*{0.2}}
\multiput(24,2)(1,0){3}{\circle*{0.2}}
\multiput(28,2)(1,0){4}{\circle*{0.2}}
\multiput(29,3)(1,0){2}{\circle*{0.2}}
\multiput(25,3)(0,2){2}{\circle*{0.2}}
\put(21,3){\circle*{0.2}}
\multiput(20,2)(4,0){2}{\line(1,0){2}}
\multiput(21,2)(4,0){2}{\line(0,1){1}}
\multiput(29,2)(1,0){2}{\line(0,1){1}}
\put(28,2){\line(1,0){3}}
\dashline{0.2}(22,2)(24,2)
\dashline{0.2}(26,2)(28,2)
\dashline{0.2}(25,3)(25,5)
\put(22,1){\small{$k-3$}}
\put(26,1){\small{$k-2$}}
\put(25.5,4.2){\small{$k-3$}}
\put(15,-.5){$(b)$}
\end{picture}
\end{center}
\caption{Open quipus with diameter $2k$}\label{ad4}
\end{figure}

By Lemma~\ref{SpD}, the only closed quipu with diameter $2(n-1)/3$ and spectral radius less than $3/\sqrt 2$ is $C_{(2k+1)}^{(k-1)}$. By Lemmas~\ref{lm2.4}~and~\ref{lm2.2}, and Corollary~\ref{trgequ2}, we get
$$\rho\left (C_{(2k+1)}^{(k-1)}\right )>\rho\left (C_{(2k+2)}^{(k-1)}\right )>\rho\left (P_{(k+1,k+1)}^{(k-1)}\right )=\rho\left (P_{(k,k)}^{(k)}\right )=\rho_k.$$
This shows that $C_{(2k+1)}^{(k-1)}$ cannot be a minimizer graph. Note that a dagger of order $n$ has diameter $n-3>2(n-1)/3$ for $n>11$. Then by Lemma~\ref{WandN}, any minimizer graph must be an open quipu with spectral radius less than $3/\sqrt 2$, which can be written as $P_{(m_0,k_1,...,k_r,m_r)}^{(m_0,m_1,...,m_r)}$. Counting the number of vertices and the diameter, we have
\begin{eqnarray}
3k&=&n-1=m_0+m_r+\sum_{j=0}^{r} m_j+\sum_{i=1}^r k_i+r,\label{n}\\
2k&=&m_0+m_r+\sum_{i=1}^r k_i+r.\label{D}
\end{eqnarray}
By Lemma~\ref{sufc}, we also have
\begin{eqnarray}
l_1 &:=& k_1 + 2-m_0 - m_1\ge 0,\label{l1}\\
l_r &:=& k_r + 2-m_{r-1} - m_r\ge 0,\mbox{ and}\label{lr}\\
l_i &:=& k_i-m_{i-1} - m_i\ge 0\mbox{ for }2\le i\le r-1.\label{li}
\end{eqnarray}
Summing up these equalities and applying Eqs.~(\ref{n}) and (\ref{D}), we obtain
\begin{eqnarray*}
\sum_{j=1}^r l_j&=&\sum_{i=1}^{r} k_i-m_0-m_r-2 \sum_{l=1}^{r-1} m_l + 4\\
&=& 3\left(m_0+m_r+\sum_{i=1}^r k_i+r\right)-2\left(m_0+m_r+\sum_{j=0}^{r} m_j+\sum_{i=1}^r k_i+r\right)+4-r\\
&=& 4-r.
\end{eqnarray*}
This implies that $r\le 4$. We will show that all open quipus with $r>1$ internal paths must have spectral radius greater than $\rho _k$, which implies the right minimizer graphs as desired in the theorem. For this purpose, those open quipus with spectral radius at most $3/\sqrt 2$ need only to be considered. One can check that Lemmas~\ref{sufc}~and~\ref{sufc2} exclude most open quipus for minimizer graphs except those shown in Figure~\ref{ad4} whose spectral radii, however, are indeed greater than $\rho _k$, as proven in the following.
\medskip

{\bf Case 1 } $r=2$. In this case, $l_1+l_2=2$. By symmetry, we have the following two subcases.
\medskip

{\bf Subcase 1.1 } $l_1=0$ and $l_2=2$.
\medskip

Eqs. (\ref{l1}) and (\ref{lr}) imply that
\begin{eqnarray*}
k_1&=&m_0+m_1-2,\\
k_2&=&m_1+m_2.
\end{eqnarray*}
Then by Lemma~\ref{sufc} (2), we have
$$\left\lceil\frac{3\delta_1(m_0)+\delta_1(m_1)}{2}\right\rceil
+\left\lfloor\frac{\delta_1(m_0-1)+\delta_1(m_1-1)}{2}\right\rfloor\ge 2.$$
It follows that $m_0=1$. Also by Lemma~\ref{sufc2},
$$\left\lceil [3\delta_1( m_2) +
\delta_1( m_1) +\delta_1( m_2-1) ]/2\right\rceil >0,$$
which implies that $m_1=1$ or $m_2=1, 2$. Combining with Eqs.~(\ref{n})~and~(\ref{D}), we obtain that all open quipus $P_{(m_0,k_1,k_2,m_2)}^{(m_0,m_1,m_2)}$, except $P_{(1,k-3,k-1,1)}^{(1,k-2,1)}$, $P_{(1,k-4,k-1,2)}^{(1,k-3,2)}$, and $P_{(1,0,k-1,k-2)}^{(1,1,k-2)}$ shown in Figure~\ref{ad4} (a), have spectral radius greater than $3/\sqrt 2$. By Lemmas~\ref{trans} and~\ref{mm'}, however, we have
\begin{eqnarray*}
\rho\left (P_{(1,k-3,k-1,1)}^{(1,k-2,1)}\right )&>&\rho\left (P_{(1,k-1,k-1)}^{(1,k-1)}\right )=\rho _k,\\
\rho\left (P_{(1,k-4,k-1,2)}^{(1,k-3,2)}\right )&>&\rho\left (P_{(2,k-1,k-2)}^{(2,k-2)}\right )=\rho _k,\\
\rho\left (P_{(1,0,k-1,k-2)}^{(1,1,k-2)}\right )&>&\rho\left (P_{(2,k-1,k-2)}^{(2,k-2)}\right )=\rho _k.
\end{eqnarray*}

{\bf Subcase 1.2 } $l_1=1$ and $l_2=1$.
\medskip

Eqs. (\ref{l1}) and (\ref{lr}) imply that
\begin{eqnarray*}
k_1&=&m_0+m_1-1,\\
k_2&=&m_1+m_2-1.
\end{eqnarray*}
Then by Lemma~\ref{sufc},
\begin{eqnarray*}
\left\lceil\frac{3\delta_1(m_0)+\delta_1(m_1)}{2} \right\rceil
+\left\lfloor\frac{\delta_1(m_0-1)+\delta_1(m_1-1)}{2}\right\rfloor &\ge & 1,\\
\left\lceil\frac{3\delta_1(m_2)+\delta_1(m_{1})}{2} \right\rceil
+\left\lfloor\frac{\delta_1(m_2-1)+\delta_1(m_1-1)}{2}\right\rfloor &\ge & 1.
\end{eqnarray*}
It follows that $m_0=m_2=1$ or $m_1=1$ since $n>20$. Also by Lemma~\ref{sufc2},
\begin{eqnarray*}
\left\lceil [3\delta_1( m_{0}) +
\delta_1( m_{1}) +\delta_1( m_{0}-1)]/2\right\rceil & > & 1,\mbox{ or}\\
\left\lceil [3\delta_1( m_{2}) +
\delta_1( m_{1}) +\delta_1( m_{2}-1)]/2\right\rceil & > & 1.
\end{eqnarray*}
It follows that $m_0=1$ or $m_2=1$. Therefore, combining with Eqs. (\ref{n}) and (\ref{D}), we obtain that all open quipus $P_{(m_0,k_1,k_2,m_2)}^{(m_0,m_1,m_2)}$, except $P_{(1,k-2,k-2,1)}^{(1,k-2,1)}$ and $P_{(1,1,k-2,k-2)}^{(1,1,k-2)}$ shown in Figure~\ref{ad4} (a), have spectral radius greater than $3/\sqrt 2$.

By Lemmas~\ref{lm2.4} and~\ref{mm'}, and Corollary~\ref{trgequ2}, we get
$$\rho\left (P_{(1,k-2,k-2,1)}^{(1,k-2,1)}\right )>\rho\left (P_{(1,k,k,1)}^{(1,k-2,1)} \right )=\rho\left (P_{(1,k-1,k-1)}^{(1,k-1)}\right )=\rho _k.$$
Let $G:=P_{(1,1,k-3)}^{(1,1)}$ and $H:=P_{(2,k-2)}^{(2)}$, and let $x$ and $y$ be the right most endvertices of $G$ and $H$ respectively. Note that
\begin{eqnarray*}
G&\cong &\left (P_{(1,3)}^{(1)},u, k-3\right ),\\
H&\cong &\left (P_{(2,1)}^{(2)},v, k-3\right ).
\end{eqnarray*}
By Lemma~\ref{compare} and Eq.~(\ref{alphaequ2}), we have $\alpha _{(G,x)}(\lambda )<\alpha _{(H,y)}(\lambda )$ for $\lambda >2$. Then by Lemmas~\ref{trg} and~\ref{mm'}, we have
$$\rho\left(P_{(1,1,k-2,k-2)}^{(1,1,k-2)}\right )>\rho\left (P_{(2,k-1,k-2)}^{(2,k-2)}\right )=\rho _k,$$
noting that
\begin{eqnarray*}
P_{(1,1,k-2,k-2)}^{(1,1,k-2)}&\cong &T_{G,P_{k-3}}^{P_{k-2}},\\
P_{(2,k-1,k-2)}^{(2,k-2)}&\cong &T_{H,P_{k-3}}^{P_{k-2}}.
\end{eqnarray*}

{\bf Case 2 } $r=3$. We have $l_1+l_2+l_3=1$, which implies that only one of $l_1$, $l_2$, and $l_3$ equals one. By symmetry, we have the following two subcases.
\medskip

{\bf Subcase 2.1 } $l_1=l_3=0$ and $l_2=1$.
\medskip

Eqs. (\ref{l1}), (\ref{lr}) and (\ref{li}) imply that
\begin{eqnarray*}
k_1&=&m_0+m_1-2,\\
k_2&=&m_1+m_2+1,\\
k_3&=&m_2+m_3-2.
\end{eqnarray*}
Then by Lemma~\ref{sufc}, we have
\begin{eqnarray*}
\left\lceil\frac{3\delta_1(m_0)+\delta_1(m_1)}{2} \right\rceil
+\left\lfloor\frac{\delta_1(m_0-1)+\delta_1(m_1-1)}{2}\right\rfloor &\ge &2,\\
\left\lceil\frac{3\delta_1(m_3)+\delta_1(m_2)}{2} \right\rceil
+\left\lfloor\frac{\delta_1(m_3-1)+\delta_1(m_2-1)}{2}\right\rfloor &\ge &2.
\end{eqnarray*}
It follows that $m_0=m_3=1$. Lemma~\ref{sufc2} however, implies that all open quipus $P_{(1,k_1,k_2,k_3,1)}^{(1,m_1,m_2,1)}$ have spectral radius greater than $3/\sqrt 2$.
\medskip

{\bf Subcase 2.2 } $l_1=1$ and $l_2=l_3=0$.
\medskip

Eqs. (\ref{l1}), (\ref{lr}) and (\ref{li}) imply that
\begin{eqnarray*}
k_1&=&m_0+m_1-1,\\
k_2&=&m_1+m_2,\\
k_3&=&m_2+m_3-2.
\end{eqnarray*}
Then by Lemma~\ref{sufc}, we have
\begin{eqnarray*}
\left\lceil [\delta_1(m_1)+\delta_1(m_2)]/2\right\rceil &\ge & 1,\\
\left\lceil\frac{3\delta_1(m_3)+\delta_1(m_2)}{2} \right\rceil
+\left\lfloor\frac{\delta_1(m_3-1)+\delta_1(m_2-1)}{2}\right\rfloor &\ge & 2.
\end{eqnarray*}
It follows that $m_1=1$ or $m_2=1$, and $m_3=1$. Also by Lemma~\ref{sufc2},
$$\left\lceil [3\delta_1( m_{0})+\delta_1( m_{1}) +\delta_1( m_{0}-1)]/2\right\rceil >1,$$
which implies that $m_0=1$. Therefore, combining with Eqs. (\ref{n}) and (\ref{D}), we obtain that all open quipus $P_{(m_0,k_1,k_2,k_3,m_3)}^{(m_0,m_1,m_2,m_3)}$, except $P_{(1,1,k-2,k-4,1)}^{(1,1,k-3,1)}$ and $P_{(1,k-3,k-2,0,1)}^{(1,k-3,1,1)}$ shown in Figure~\ref{ad4} (b), have spectral radius greater than $3/\sqrt 2$.

By Lemma~\ref{trans}, however, we get
$$\rho\left (P_{(1,1,k-2,k-4,1)}^{(1,1,k-3,1)}\right )>\rho\left (P_{(1,1,k-2,k-2)}^{(1,1,k-2)}\right )>\rho _k,$$ where the last inequality was proved in Subcase 1.2, and
$$\rho\left (P_{(1,k-3,k-2,0,1)}^{(1,k-3,1,1)}\right )
>\rho\left (P_{(1,k-3,k-2,2)}^{(1,k-3,2)}\right )
>3/\sqrt 2,$$
where the last inequality holds since $k-2<k-3+2$ and $\delta _1(k-3)=\delta _1(k-4)=0$ for $k>6$ which fails to satisfy Lemma~\ref{sufc} (3).
\medskip

{\bf Case 3 } $r=4$. We have $l_1+l_2+l_3+l_4=0$, which implies that $l_1=l_2=l_3=l_4=0$. Eqs. (\ref{l1}), (\ref{lr}) and (\ref{li}) imply that
\begin{eqnarray*}
k_1&=&m_0+m_1-2,\\
k_2&=&m_1+m_2,\\
k_3&=&m_2+m_3,\\
k_4&=&m_3+m_4-2.
\end{eqnarray*}
As above, Lemma~\ref{sufc} implies that $m_0=m_4=1$. Lemma~\ref{sufc2} however, implies that all open quipus $P_{(1,k_1,k_2,k_3,k_4,1)}^{(1,m_1,m_2,m_31)}$ have spectral radius greater than $3/\sqrt 2$. This completes the proof.
$\hfill\Box $

\end{document}